\documentclass[11pt]{article}

\usepackage{amsfonts,amsmath,amssymb,cite}
\usepackage{amsmath,amsthm}
\usepackage{enumerate}
\usepackage{latexsym}
\usepackage{color}
\usepackage[dvipsnames]{xcolor}
\usepackage{tikz}
\usetikzlibrary{decorations.pathmorphing}
\usetikzlibrary{shapes.geometric}
\usetikzlibrary{svg.path}
\usetikzlibrary{arrows,positioning,backgrounds}

\newtheorem{thm}{Theorem}

\newtheorem{ob}{Observation}

\newtheorem{cor}[thm]{Corollary}
\newtheorem{prop}[thm]{Proposition}

\theoremstyle{definition}

\theoremstyle{remark}

\newcommand{\ggr}{\gamma_{\rm gr}}

\newcommand{\vertex}{\node[vertex]}
\tikzstyle{vertex}=[circle, draw, inner sep=0pt, minimum size=6pt]

\begin{document}

\title{Loop zero forcing and grundy domination in planar graphs and claw-free cubic graphs}
\author{
Alex Domat$^{a}$ \and
Kirsti Kuenzel$^{a}$ }

\maketitle

\begin{center}
$^a$ Department of Mathematics, Trinity College, Hartford, CT, USA\\
\end{center}
\medskip

\maketitle
\begin{abstract}
Given a simple, finite graph with vertex set $V(G)$, we define a zero forcing set of $G$ as follows.  Choose $S\subseteq V(G)$ and  color all vertices of $S$ blue and all vertices in $V(G) - S$ white. The color change rule is if $w$ is the only white neighbor of blue vertex $v$, then we change the color of $w$ from white to blue. If after applying the color change rule as many times as possible eventually every vertex of $G$ is blue, we call $S$ a zero forcing set of $G$. $Z(G)$ denotes the minimum cardinality of a zero forcing set. Davila and Henning proved in \cite{zerocubic} that for any claw-free cubic graph $G$, $Z(G) \le \frac{1}{3}|V(G)| + 1$. We show that if $G$ is $2$-edge-connected, claw-free, and cubic, then $Z(G) \le \left\lceil\frac{5n(G)}{18}\right\rceil+1$. We also study a similar graph invariant known as the loop zero forcing number of a graph $G$ which happens to be the dual invariant to the Grundy domination number of $G$. Specifically, we study the loop zero forcing number in two particular types of planar graphs. 

\end{abstract}

%%%%%%%%%%%%%%%%%%%%%%%%%%%%%%%%%%%
%%%%%%%%%%%%%%%%%%%%%%%%%%%%%%%%%%%
\section{Introduction} \label{sec:intro}
Zero forcing is a coloring process of the vertices of $G$ according to the following rule. Pick a set $S \subseteq V(G)$ and color each vertex of $S$ blue and each vertex of $V(G) - S$ white. The color change rule is that for each $w \in V(G) - S$ where $w$ is the only white neighbor of a blue vertex $v$, then color $w$ blue. If after applying the color change rule as many times as possible each vertex of $G$ is blue, then we call $S$ a zero forcing set of $G$. The zero forcing number of $G$, denoted $Z(G)$, is the minimum cardinality among all zero forcing sets of $G$. The original motivation for studying zero forcing was that $Z(G)$ is an upper bound for the minimum rank of a graph (see \cite{AIM}). In a seemingly unrelated paper, Bre\v{s}ar et al. \cite{domseq} defined what is known as a \emph{legal dominating sequence} and \emph{Grundy domination number} of a graph in 2014  as follows. Let $S = (x_1, \dots, x_n)$ be a sequence of distinct vertices of $G$. $S$ is called a legal sequence if for each $i \in [n-1]$, $x_i$ dominates some vertex that is not already dominated by any vertex in $\{x_1, \dots, x_{i-1}\}$. If $S= (x_1, \dots, x_n)$ is a legal sequence and $\{x_1, \dots, x_n\}$ dominates $G$, then we call $S$ a legal dominating sequence of $G$. A longest possible dominating sequence is called a Grundy dominating sequence and its length is denoted $\ggr(G)$. Bre\v{s}ar et al. showed in  \cite{Zgrundyzero} that by modifying this definition slightly, one derives the dual invariant of $Z(G)$. In particular, $S = (x_1, \dots, x_n)$ is said to be a legal $Z$-sequence if for each $i \in [n-1]$, $x_i$ has a neighbor that is not dominated by any vertex in $\{x_1, \dots, x_{i-1}\}$. The longest $Z$-sequence in $G$ is called the $Z$-Grundy domination number and its length is denoted $\ggr^Z(G)$. As $\ggr^Z(G)$ is the dual invariant of $Z(G)$, it follows that  $|V(G)| = \ggr^Z(G) + Z(G)$. 

Upon realizing that a Grundy-like domination number is the dual invariant to the zero forcing number of a graph, Lin \cite{grundyloopzero} showed that the Grundy domination number of a graph is the dual invariant of what is called the loop zero forcing number of a graph, defined as follows. In \cite{minrank}, the zero forcing number of a loop graph, denoted $Z_{\dot{\ell}}(G)$,  was defined whereby a loop graph is constructed from any simple graph $G$ by adding a loop at each vertex of $G$. On the other hand, the loop zero forcing number of a graph, $Z_{\ell}(G)$,  was defined in \cite{loopzero} slightly differently as finding the zero forcing number of the graph $\widehat{G}$ which is obtained from a simple graph $G$ by adding a loop at each vertex $v \in V(G)$ where $\deg_G(v) \ge 1$. Note however, that if we only consider connected, simple graphs $G$, $Z_{\ell}(G) = Z_{\dot{\ell}}(G)$. Therefore, we refer to $Z_{\dot{\ell}}(G)$ as the loop zero forcing number of $G$ as we only consider connected, simple graphs throughout this paper. Moreover, for any connected, simple graph $G$, $|V(G)| = \ggr(G) + Z_{\dot{\ell}}(G)$. 

Davila and Henning showed in \cite{zerocubic} that for any claw-free cubic graph $G$ of order at least $10$ that $Z(G) \le \frac{1}{3}|V(G)| +1$. Motivated by their result, we consider the loop zero forcing number of $2$-connected, claw-free cubic graphs. We use a similar technique used by Anderson and Kuenzel in \cite{pdcubic} whereby  an upper bound for the power domination number of a claw-free cubic graph was given. Although originally motivated by loop zero forcing sets, we show the following upper bound for $Z(G)$ in the case where $G$ is a $2$-connected, claw-free cubic graph.

\begin{thm}\label{thm:cubic} Let $G$ be a $2$-edge-connected, claw-free cubic graph. If $G$ is not a ring of diamonds, then $Z_{\dot{\ell}}(G)\le Z(G) \le \left\lceil\frac{5n(G)}{18}\right\rceil+1$. If $G$ is a ring of diamonds, then $Z_{\dot{\ell}}(G) =Z(G) =  \frac{n(G)}{4}+2$. 
\end{thm}

We also consider the loop zero forcing number of certain types of planar graphs. In \cite{triangulations}, a series of upper bounds for $Z(G)$ were provided when $G$ is a maximal outerplanar graph. Motivated by their work, we show that in certain cases the loop zero forcing number and zero forcing number of maximal outerplanar graphs coincide. We also provide a lower bound for the zero forcing number of a maximal outerplanar graph. 

The remainder of this paper is organized as follows. In Section~\ref{sec:def}, we provide definitions and terminology used throughout the paper. Section~\ref{sec:planar} is dedicated to studying the loop zero forcing number of certain planar graphs. In particular, we consider maximal outerplanar graphs and Halin graphs. In Section~\ref{sec:cubic}, we prove Theorem~\ref{thm:cubic}. 
%%%%%%%%%%%%%%%%%%%%%%%%%%%%%%%%%%%
%%%%%%%%%%%%%%%%%%%%%%%%%%%%%%%%%%%

\subsection{Definitions and Preliminaries}\label{sec:def}
Let $G$ be a finite, simple graph with vertex set $V(G)$ and edge set $E(G)$.  We let $n(G) = |V(G)|$.  For a vertex $x \in V$, the {\em open neighborhood of $x$}
is the set $N(x)$ defined by $N(x)=\{w \in V(G):\, xw \in E(G)\}$.  The {\em closed neighborhood} $N[x]$ is
$N(x) \cup \{x\}$.   The open neighborhood of a set $A \subseteq V(G)$ is $N(A)=\cup_{a\in A}N(a)$ and its
closed neighborhood is $N[A]=N(A) \cup A$. Given two vertices $x, y\in V(G)$, we let $d_G(x, y)$ denote the distance between $x$ and $y$, or length of the shortest path between $x$ and $y$. 

A set $A$ of vertices  is a {\em dominating} set of $G$ if $N[A]=V(G)$.  Let $S=(x_1,\ldots,x_n)$ be a sequence of distinct vertices in $G$.  We denote the length of $S$ by $|S|$.  The set $\{x_1,\ldots,x_n\}$ whose elements are the vertices in $S$ is denoted  by $\widehat{S}$. The sequence $S$ is called a  {\em closed neighborhood sequence} (or a {\em legal sequence}) if
\begin{equation}
\label{eq:defGrundy}
N[x_{i+1}]-\bigcup_{j=1}^{i}N[x_j]\not=\emptyset
\end{equation}
for each $i\in [n-1]$.  We will also say that $x_{i+1}$ {\em footprints} the vertices from $N[x_{i+1}] - \bigcup_{j=1}^{i}N[x_j]$ with respect to $S$.   If $S$ is a legal sequence and $\widehat{S}$ is a dominating set of $G$, then $S$ is called a {\em dominating sequence} in $G$.  A longest possible dominating sequence in $G$ is called a {\em Grundy dominating sequence},  and its length is the {\em Grundy domination number} of $G$, denoted $\ggr(G)$.

A closed neighborhood sequence $S$ in $G$ is a $Z$-\emph{sequence} if, in addition, every vertex $v_i$ in $S$ footprints a vertex distinct from itself. That is, $S$ is a Z-sequence in $G$ if 
\[N_G(v_i) \setminus \bigcup_{j=1}^{i-1}N_G[v_j] \ne \emptyset\]
holds for every $i \in \{2, \dots, k\}$. The maximum length of a Z-sequence in $G$ is the \emph{Z-Grundy domination number}, denoted $\gamma_{\rm{gr}}^Z(G)$, of $G$. Given a Z-sequence $S$, we call the corresponding set $\widehat{S}$ a Z-set. 

Given a set $S\subseteq V(G)$, we define $S$ to be a zero forcing set of $G$ if the following is true. Initially, color all vertices of $S$ blue and all vertices of $V(G) - S$ white. We iteratively apply the zero forcing color change rule: if $w \in V(G)$ is the only white neighbor of blue vertex $v$, recolor $w$ blue. If after applying the zero forcing color change rule as many times as possible all vertices of $G$ are blue, we call $S$ a \emph{zero forcing set} of $G$. $Z(G)$ denotes the minimum cardinality among all zero forcing sets of $G$.  We will say $S\subseteq V(G)$ is a loop zero forcing set if the following is true. Initially, we color all vertices of $S$ blue and all vertices of $V(G) - S$ white. We iteratively apply the loop zero forcing color change rule: if $w\in V(G)$ is either 1) the only white neighbor of blue vertex $v$, or 2) every neighbor of $w$ is blue, then we recolor $w$ blue. If after applying the loop zero forcing color change rule as many times as possible all vertices of $G$ are blue, then call $S$ a \emph{loop zero forcing set} of $G$. $Z_{\dot{\ell}}(G)$ denotes the minimum cardinality among all loop zero forcing sets of $G$. Note that the zero forcing number of $G$ is the dual invariant of $\ggr^Z(G)$ and the loop zero forcing number of $G$ is the dual invariant of $\ggr(G)$. That is, $n(G) = Z(G) + \ggr^Z(G) = Z_{\dot{\ell}}(G) +\ggr(G)$. 

We will make use of the following two known results. 
\begin{prop}\label{prop:mindegree}\cite{domseq} For an arbitrary graph $G$, $\ggr(G) \le n(G) - \delta(G)$. 
\end{prop}

\begin{thm}\label{thm:bresar}\cite{Zgrundyzero} If $G$ is a graph without isolated vertices, then $\ggr^Z(G) + Z(G) = n(G)$. Moreover, the complement of a zero forcing set of $G$ is a $Z$-set of $G$ and vice versa. 
\end{thm}

%%%%%%%%%%%%%%%%%%%%%%%%%%%%%%%%%%%%%%%%%%%%%%%%%%%%%%%%%%%%%%%%%%%%%%%%%%%%%%%%%%%%%%%%%%%%%%%%%%%%%%%%%%%%%%%%%%%%%%%%%%%%%%%%%%%%%%%%%%%%%%%%

\section{Planar Graphs}\label{sec:planar}
In this section, we consider two types of planar graphs. Recall that $G$ is \emph{outerplanar} if $G$ has a planar drawing in which all vertices of $V(G)$ belong to the outer face of the drawing. The \emph{weak planar dual} graph of a planar graph is the graph that has a vertex for every bounded face and an edge between two vertices which represent a pair of adjacent bounded faces. In this section, we consider both maximal outerplanar graphs and Halin graphs.

\subsection{Maximal Outerplanar Graphs}
\begin{figure}[h]
\begin{center}
\begin{tikzpicture}[scale=1]
\tikzstyle{vertex}=[circle, draw, inner sep=0pt, minimum size=6pt]
\tikzset{vertexStyle/.append style={rectangle}}
	\vertex (1) at (0,-.5) [scale=.75, label=below:$$] {};
	\vertex (2) at (-.5, .25) [scale=.75, label=below:$$] {};
	\vertex (3) at (-.5, 1) [scale=.75, label=below:$$] {};
	\vertex (4) at (-.5, 1.75) [scale=.75, label=below:$$] {};
	\vertex (5) at (-1.5, 1.75) [scale=.75, label=below:$$] {};
	\vertex (6) at (-2.5, 1.75) [scale=.75, label=below:$$] {};
	\vertex (7) at (-3, 2.5) [scale=.75, label=below:$$] {};
	\vertex (8) at (-2, 2.5) [scale=.75, label=below:$$] {};
	\vertex (9) at (-1, 2.5) [scale=.75, label=below:$$] {};
	\vertex (10) at (0,2.5) [scale=.75, label=below:$$] {};
	\vertex (11) at (1, 2.5) [scale=.75, label=below:$$] {};
	\vertex (12) at (2, 2.5) [scale=.75, label=below:$$] {};
	\vertex (13) at (3, 2.5) [scale=.75, label=below:$$] {};
	\vertex (14) at (4, 2.5) [scale=.75, label=below:$$] {};
	\vertex (15) at (5, 2.5) [scale=.75, label=below:$$] {};
	\vertex (16) at (4.5, 1.75) [scale=.75, label=below:$$] {};
	\vertex (17) at (3.5, 1.75) [scale=.75, label=below:$$] {};
	\vertex (18) at (3.5, 1) [scale=.75, label=below:$$] {};
	\vertex (19) at (3.5, .25) [scale=.75, label=below:$$] {};
	\vertex (20) at (3, -.5) [scale=.75, label=below:$$] {};
	\vertex (21) at (2.5, .25) [scale=.75, label=below:$$] {};
	\vertex (22) at (2.5, 1) [scale=.75, label=below:$$] {};
	\vertex (23) at (2.5, 1.75) [scale=.75, label=below:$$] {};
	\vertex (24) at (1.5, 1.75) [scale=.75, label=below:$$] {};
	\vertex (25) at (.5, 1.75) [scale=.75, label=below:$$] {};
	\vertex (26) at (.5, 1) [scale=.75, label=below:$$] {};
	\vertex (27) at (.5, .25) [scale=.75, label=below:$$] {};
	
	\begin{scope}[on background layer]    
 \fill[fill=green!40]   (-.5,1.75)  to (0,2.5) to (.5, 1.75) to  (-.5, 1.75)  ;
 \fill[fill=green!40]   (2.5, 1.75)  to (3, 2.5) to (3.5, 1.75) to  (2.5, 1.75)  ;
  \fill[fill=purple!40]   (3, 2.5)  to (4,2.5) to (3.5, 1.75) to  (3, 2.5)  ;
   \fill[fill=purple!40]   (3.5,1.75)  to (4, 2.5) to (4.5, 1.75) to  (3.5, 1.75)  ;
    \fill[fill=purple!40]   (4, 2.5)  to (5, 2.5) to (4.5, 1.75) to  (4, 2.5)  ;
     \fill[fill=yellow!40]   (0,2.5)  to (1,2.5) to (.5, 1.75) to  (0, 2.5)  ;
      \fill[fill=yellow!40]   (.5, 1.75)  to (1, 2.5) to (1.5, 1.75) to  (.5, 1.75)  ;
       \fill[fill=yellow!40]   (1, 2.5)  to (2, 2.5) to (1.5, 1.75) to  (1, 2.5)  ;
        \fill[fill=yellow!40]   (1.5, 1.75)  to (2, 2.5) to (2.5, 1.75) to  (1.5, 1.75)  ;
         \fill[fill=yellow!40]   (2, 2.5)  to (3, 2.5) to (2.5, 1.75) to  (2, 2.5)  ;

\end{scope}

	\path

	(1) edge (2)
	(2) edge (3)
	(3) edge (4)
	(4) edge (5)
	(5) edge (6)
	(6) edge (7)
	(7) edge (8)
	(8) edge (9)
	(9) edge (10)
	(10) edge (11)
	(11) edge (12)
	(12) edge (13)
	(13) edge (14)
	(14) edge (15)
	(15) edge (16)
	(16) edge (17)
	(17) edge (18)
	(18) edge (19)
	(19) edge (20)
	(20) edge (21)
	(21) edge (22)
	(22) edge (23)
	(23) edge (24)
	(24) edge (25)
	(25) edge (26)
	(26) edge (27)
	(1) edge (27)
	(2) edge (27)
	(2) edge (26)
	(2) edge (25)
	(3) edge (25)
	(4) edge (25)
	(4) edge (10)
	(4) edge (9)
	(5) edge (9)
	(5) edge (8)
	(6) edge (8)
	(25) edge (11)
	(11) edge (24)
	(12) edge (24)
	(12) edge (23)
	(23) edge (13)
	(13) edge (17)
	(17) edge (14)
	(14) edge (16)
	(21) edge (19)
	(21) edge (18)
	(18) edge (22)
	(22) edge (17)
	(17) edge (23)
	(10) edge (25)

	;
\end{tikzpicture}
\end{center}
\caption{Example of serpentine path and fan leaf in a MOP}
\label{fig:notation}
\end{figure}
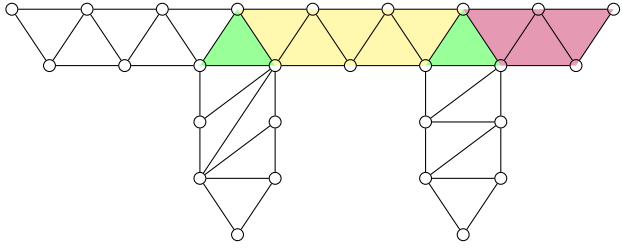

$G$ is \emph{maximal outerplanar}, or a MOP, if $G$ is outerplanar and the addition of any edge results in a graph that is not outerplanar. Thus, every interior face of $G$ is a triangle. We use similar terminology to that found in \cite{triangulations}. Namely, a triangular face $T$ is a separator triangle of $G$ if it has no edges on the outer face. The number of separator triangles of $G$ is denoted by $t$. We say a MOP $G$ is \emph{serpentine} if it contains no separator triangles. We let $G_{\triangle}$ denote the subgraph formed by the separator triangles of $G$. The weak planar dual of a MOP $G$ is a tree that we denote by $H$. A \emph{serpentine leaf} is a set of triangles in $G$ whose corresponding vertices in $H$ form the only path in $H$ connecting a leaf of $H$ with a vertex $u$ of $H$ where $\deg_H(u) >2$. We let $h$ represent the number of serpentine leaves in $G$ and $n_2$ represent the number of vertices of degree $2$ in $G$. Note that $h= n_2$. Similarly, we define a \emph{serpentine path} to be a maximal set of triangles in $G$ whose corresponding vertices in $H$ have degree $2$ in $H$ and form the only path between two distinct vertices of degree $3$ or more in $H$. The number of serprentine leaves with fan structure from a vertex of a separator triangle, called fan leaves, is denoted by $h_F$. Figure~\ref{fig:notation} provides an example of a MOP where the green triangles are separator triangles, the purple triangles lie form a fan leaf, and the yellow triangles form a serpentine path. 
\vskip5mm

Note that the following were shown in \cite{triangulations}.

\begin{prop}\cite{triangulations}\label{prop:serp} If $G$ is a serpentine graph, then $Z(G) = 2$. 
\end{prop}

\begin{thm}\cite{triangulations}\label{thm:one} If $G$ is a MOP with $t=1$ and $h_F\ge 0$, then 
\begin{enumerate}
\item[1.] If $h_F \ge 1$, then $Z(G) = 3$. 
\item[2.] If $h_F=0$, then $Z(G) = 4$. 
\end{enumerate}
\end{thm}

\begin{thm}\cite{triangulations}\label{thm:upper} If $G$ is a MOP with $n_2$ vertices of degree $2$ whose $G_{\triangle}$ has $c$ components, then 
\[Z(G) \le 2\left\lceil\frac{n_2+c-1}{2}\right\rceil\]
and these bounds are tight when the number of fan leafs is $0$ and if $c=0$ or $G_{\triangle}$ is formed by isolated triangles. 
\end{thm}

We first point out some immediate consequences of the above results. First, note that if $G$ is a serpentine graph, then $\delta(G)=2$. It is clear that for any graph $H$, $Z_{\dot{\ell}}(H)\ge \delta(H)$ and therefore by Proposition~\ref{prop:serp} we have the following. 
\begin{cor}\label{cor:serp} If $G$ is a serpentine graph, then $Z_{\dot{\ell}}(G) =2$. 
\end{cor}

Next, we show that for all MOPs $G$ described in the statement of Theorem~\ref{thm:one}, $Z_{\dot{\ell}}(G) = Z(G)$. 

\begin{prop} If $G$ is a MOP with $t=1$ and $h_F \ge 0$, then 
\begin{enumerate}
\item[1.] If $h_F \ge 1$, then $Z_{\dot{\ell}}(G) =3$.
\item[2.] If $h_F=0$, then $Z_{\dot{\ell}}(G) =4$.
\end{enumerate}
\end{prop}
\begin{proof}
First, assume that $G$ is a MOP with $t=1$ and $h_F\ge 1$. It is clear that $Z_{\dot{\ell}}(G) \ge 2$ and every vertex of $G$ is contained on a serpentine leaf. Suppose that there exists a loop zero forcing set of cardinality $2$, say $\{x, y\}$. Thus, both $x$ and $y$ must be vertices on the same serpentine leaf, call it $L$. Let $w$ and $z$ be the two vertices of $G$ that are on the one separator triangle of $G$ and also contained in $L$. Assuming that all vertices of $L$ eventually turn blue, each of $w$ and $z$ have at least two neighbors that are not in $L$ and therefore $\{x, y\}$ is not a loop zero forcing set of $G$. 

Next, assume that $G$ is a MOP with $t=1$ and $h_F = 0$. Enumerate the serpentine leaves of $G$ as $L_1, L_2, L_3$. Suppose that there exists a loop zero forcing set $S = \{x, y, t\}$. Thus, at least two vertices of $S$ must be contained in the same serpentine leaf. Without loss of generality, we may assume $x$ and $y$ are vertices in $L_1$. As above, let $w$ and $z$ be the two vertices of $G$ that are on the one separator triangle of $G$ and also contained in $L_1$. Assuming that all vertices of $L_1$ eventually turn blue, each of $w$ and $z$ have two neighbors that are not in $L_1$. Therefore, $t$ is a neighbor of at least one of $w$ or $z$. Assume first that $t$ is on the separator triangle. We shall assume $w$ and $t$ are vertices of $L_2$ and $z$ and $t$ are vertices of $L_3$. If $w$ has only one white neighbor, say $a$,  in $L_2$ and $z$ has only one white neighbor in $L_3$, say $b$, then $a$ and $b$ both turn blue and the propagation process stops as each of $a$, $b$, and $t$ have at least two white neighbors. So we shall assume that $w$ has two white neighbors in $L_2$, say $a$ and $c$. If $z$ has only one white neighbor, say $b$, in $L_3$, then $b$ turns blue. Moreover, no vertices in $L_2$ will turn blue until all vertices of $L_3$ are blue. This will only occur if $L_3$ is a fan leaf, which contradicts the assumption that $h_F=0$. Thus, this case cannot occur and we shall assume that $t$ is not on the separator triangle and $t$ is a vertex in $L_3$ adjacent to $z$. Let $r$ be the third vertex of the separator triangle. If $rt \in E(G)$, then $r$ is the only white neighbor of $z$ and therefore $r$ turns blue. Using a similar argument as that above, no vertex of $L_2$ will turn blue until all vertices of $L_3$ are blue which can only occur if $L_3$ is a fan leaf, which is a contradiction. Therefore, we may assume $rt \not\in E(G)$ and the only blue vertices that are not contained in $V(L_1) - \{w, z\}$ are $w$, $z$, and $t$. Note that $w$ and $z$ each have at least two white neighbors. Therefore, if $S$ is indeed a loop zero forcing set, it must be that $t$ has only one white neighbor. This implies that $\deg_G(t) = 2$ and $L_3$ is a fan leaf, another contradiction. Hence, no such loop zero forcing set exists. 
\end{proof}

Next, we give a lower bound for the zero forcing number of a MOP.  We begin with the following observation.

\begin{ob} For any serpentine graph $G$, there exist only two distinct minimum zero forcing sets and only two distinct maximum $Z$-sets. 
\end{ob}

The above follows from the fact that any minimum zero forcing set contains exactly two vertices, one of which is a vertex of degree $2$ in $G$, call it $v$,  and the other which is a neighbor of $v$. Therefore, by Theorem~\ref{thm:bresar}, the only maximum $Z$-sets of $G$ are the complements of the two minimum zero forcing sets. 

\begin{thm}\label{thm:half} Let $G$ be a MOP with $t\ge 1$ where $G_{\triangle}$ contains only one component. Then $Z(G) \ge \left\lceil\frac{n_2}{2}\right\rceil$. 
\end{thm}

\begin{proof} Let $H$ be the graph obtained from $G$ by removing all vertices of $G$ that are contained in a serpentine leaf, but not contained on a separator triangle. Note that $H$ is necessarily a MOP and $H$ is a serpentine graph.  If $H$ is a triangle, then we know $Z(G) \ge 3$ by Theorem~\ref{thm:one}. Therefore, we may assume that $t\ge 2$. We enumerate the vertices of $H$ as $x_1\dots x_r$ such that $x_ix_{i+1} \in E(H)$ for $i \in [r]$. Let $S = (s_1, \dots, s_n)$ be a $Z$-Grundy sequence of $G$ and suppose $L_i$ is a serpentine leaf of $G$ such that $V(L_i) \cap \widehat{S} = V(L_i)$. Moreover, assume $x_i$ and $x_{i+1}$ are in $L_i$. Suppose first that $x_{i+1}$ and $x_{i+2}$ are also contained in a serpentine leaf of $G$, call it $L_{i+1}$. Since all vertices of $L_i$ are in $\widehat{S}$, some vertex $y \in V(L_i) - \{x_i, x_{i+1}\}$ must footprint $x_i$ and some vertex $z \in V(L_i) - \{x_i, x_{i+1}, y\}$ must footprint $x_{i+1}$ since there is only one distinct $Z$-Grundy sequence of the serpentine graph $L_i$ not containing $x_i$ and $x_{i+1}$. Now since $z$ footprints $x_{i+1}$ and $x_{i+1}\in \widehat{S}$, there exists a neighbor $t$ of $x_{i+1}$ in $L_{i+1}$ such that $t \not\in \{x_{i+1}, x_{i+2}\}$ and $t$ does not precede $x_{i+1}$ in $S$. Let $\overline{S}$ be the subsequence of $S$ that contains only those vertices in $L_{i+1}$. If $\overline{S}$ is a $Z$-Grundy sequence of the serpentine graph $L_{i+1}$, then one of $x_{i+1}$ or $x_{i+2}$ has degree $2$ in $L_{i+1}$ and the other vertex of $L_{i+1}$ that has degree $2$ is not in $\widehat{S}$ as there is only one distinct $Z$-Grundy sequence of the serpentine graph $L_{i+1}$ where $x_{i+1}$ precedes $t$. On the other hand, if $\overline{S}$ is not a $Z$-Grundy sequence of the serpentine graph $L_{i+1}$, then some vertex of $V(L_{i+1}) - \{x_{i+1}\}$ is not in $\widehat{S}$. In either case, some vertex of $L_{i+1}$ is not in $\widehat{S}$.

Therefore, we shall assume that $x_{i+1}x_{i+2}$ is on the boundary of $G$. We may extend $x_{i+1}x_{i+2}$ to a maximal path $x_{i+1}x_{i+2}\dots x_j$ that lies on the boundary of both $H$ and $G$ and where $x_j$ is contained in a serpentine leaf yet $x_{i+2}, \dots, x_{j-1}$ are not. Let $L_j$ be the serpentine leaf containing $x_j$ and assume that $V(L_j) \cap \widehat{S} = V(L_j)$. Moreover, assume that $\{x_{i+1}, \dots, x_j\} \subset \widehat{S}$. Suppose first that $x_j = x_{i+2}$. Since every vertex of $L_i$ is in $\widehat{S}$, $x_{i+2}$ comes after $x_{i+1}$ in $S$ for otherwise some vertex of $L_i$ does not footprint a vertex. However, using this same logic, since every vertex of $L_j$ is in $\widehat{S}$, $x_{i+1}$ comes after $x_{i+2}$ in $S$, which is a contradiction. Therefore, we shall assume that $j>i+2$. As above, since every vertex of $L_i$ is in $\widehat{S}$, $x_{i+2}$ comes after $x_{i+1}$ in $S$ for otherwise some vertex of $L_i$ does not footprint a vertex. Similarly, $x_{j-1}$ comes after $x_j$ in $S$. Thus, of the set $\{x_{i+1}, \dots, x_j\}$, the last vertex to appear in $S$ is in $\{x_{i+2}, \dots, x_{j-1}\}$, call it $x_{\alpha}$. But this cannot be as $x_{\alpha-1}$ and $x_{\alpha+1}$ precede $x_{\alpha}$ in $S$, meaning $x_{\alpha}$ does not footprint a vertex. Therefore, it must be that either some vertex of $\{x_{i+2}, \dots, x_{j-1}\}$ is not in $\widehat{S}$, or some vertex of $L_j$ is not in $\widehat{S}$. In either case, we have $|\widehat{S}\cap V(G)| \le n(G) - \left\lceil\frac{n_2}{2}\right\rceil$ which implies $Z(G) \ge n(G) - (n(G) -  \left\lceil\frac{n_2}{2}\right\rceil) = \left\lceil\frac{n_2}{2}\right\rceil$.

\end{proof}

We use the above result to provide a lower bound for the zero forcing number of any MOP. 

\begin{cor} If $G$ is a MOP where $t\ge 1$, $\mathcal{H}$ represents the set of components in $G_{\triangle}$ that are adjacent to exactly one serpentine path in $G$ and $n_2'$ is the number of serpentine leaves that are adjacent to a component in $\mathcal{H}$, then $Z(G) \ge \left\lfloor \frac{n_2'}{2}\right\rfloor - 2c'$ where $c'$ is the number of components in $\mathcal{H}$.
\end{cor}

\begin{proof} Let $\overline{G}$ be the graph obtained from $G$ by removing all vertices of $G$ that are contained in a serpentine leaf or path, but not contained  on a separator triangle. Let $H \in \mathcal{H}$ and let $\overline{H}$ be the component of $\overline{G}$ where $V(H) \subseteq V(\overline{H})$. Thus, $\overline{H}$ is a MOP. We may enumerate the vertices of $\overline{H}$ as $x_1\dots x_r$ where $x_ix_{i+1}$ is on the boundary of $\overline{H}$ for $i \in [r]$. Reindexing if necessary, we may assume $x_1x_2$ is on a serpentine path of $G$. Let $n_H$ be the number of serpentine leaves adjacent to $H$. Let $G'$ be the graph containing $V(H)$ obtained from $G$ by removing the neighbors of $x_1$ and $x_2$ that lie on the serpentine path adjacent to $H$ but are not themselves contained in a separator triangle. Thus, $G'$ is a MOP where $G'_{\triangle}$ contains one component. Let $S$ be a zero forcing set of $G$. Even if $x_1$ and $x_2$ are observed by vertices in their serpentine path, from Theorem~\ref{thm:half} we know $|S\cap V(G')| \ge  \left\lceil\frac{n_H}{2}\right\rceil$ which implies $|S \cap (V(G') - \{x_1, x_2\})| \ge \left\lceil\frac{n_H}{2}\right\rceil - 2$. Summing over all components in $\mathcal{H}$ yields the desired result. 

\end{proof}

\subsection{Halin Graphs}

\noindent
 Recall that a Halin graph is constructed from a tree by enumerating the leaves of $T$ as $\ell_1, \dots, \ell_m$ in such a way so that when we add the edges $\ell_i\ell_{i+1}$ for each $i \in [m]$ to $T$, the resulting graph is planar. In this section, we use a similar strategy (and therefore similar terminology) to that used by Bre\v{s}ar et al. in \cite{domseq} where the authors are providing bounds on the Grundy domination number of a tree. For any tree $T$, a support vertex of $T$ is any vertex adjacent to a leaf in $T$. An \emph{end support vertex} is a support vertex $v$  such that $\deg'(v) \leq 1$, where $\deg'(v)$ is the number of vertices adjacent to $v$ which are not leaves. Equivalently, a support vertex $v$ is an end support vertex if and only if $v$ does not lie on a path between two other support vertices. The following lower bound was given for the Grundy domination number of a tree in \cite{domseq}.
 
 \begin{prop}\cite{domseq}\label{prop:trees} For any tree $T$ which is not a star, $\ggr(T) \ge n(T) - |ES(T)| + 1$ where $ES(T)$ is the set of end support vertices in $T$.
 \end{prop}
 
 This immediately implies the following in terms of loop zero forcing. 
 
 \begin{cor}\cite{domseq} For any tree $T$ which is not a star, $Z_{\dot{\ell}}(T) \le |ES(T)| -1$. 
 \end{cor}
 
 Studying the loop zero forcing number of Halin graphs is quite different from studying the loop zero forcing number of trees as there are no leaves in a Halin graph. Given a Halin graph $G$ constructed from a tree $T_G$, we say that a vertex of $G$ is an end support vertex of $G$ if it is an end support vertex of the tree $T_G$.

\begin{prop} If $ES(G)$ is the set of end support vertices in a Halin graph $G$, then  $Z_{\dot{\ell}}(G)\le Z(G) \leq 3(|ES(G)| - 1)$.
\end{prop}

\begin{proof} Assume that $G$ is constructed from $T_G$ and let $C$ represent the exterior boundary of $G$ consisting of all the leaves in $T_G$. Let $L$ be the set of leaves of $T_G$, and let $H = G - L$. Note that $H$ is a tree. We then construct a tree $T$ with vertices $v_1, \dots, v_t$ such that each $v_i \in V(T)$ corresponds to a unique vertex of $H$ which does not have degree $2$ in $H$. $v_i$ and $v_j$ are adjacent in $T$ if their corresponding vertices in $H$ are adjacent, or the path $P$ in $H$ between the corresponding vertices in $H$ is such that all interior vertices of $P$ have degree $2$ in $H$.  Root $T$ at a leaf $r \in T$. Note that the leaves of $T$ correspond to the end support vertices of $G$. We refer to  vertices in $G$ that correspond to non-leaf vertices in $T$ as \textit{branching vertices}.

We define two kinds of branch structures in $G$. These are analogs of serpentine leaves and serpentine paths for MOPs. A type-1 branch $B_1 \subset G$ is constructed as follows. Let $xy \in E(T)$ where $x$ is a leaf of $T$ and let $u_x$ (resp. $u_y$) be the vertex of $G$ that $x$ represents (resp. $y$ represents). The vertices of $B_1$ consist of $u_x$, $u_y$, any vertices along the $u_x$-$u_y$ path in $T_G$, call it $P = u_xw_1\cdots w_ku_y$, as well as any leaves (in $T_G$) of any vertex of $P$ other than $u_y$.  A type-2 branch $B_2 \subset G$ is constructed as follows. Let $xy \in E(T)$ where neither $x$ nor $y$ is a leaf of $T$ and let $u_x$ (resp. $u_y$) be the vertex of $G$ that $x$ represents (resp. $y$ represents). The vertices of $B_2$ consist of $u_x$, $u_y$, any vertices along the $u_x$-$u_y$ path in $T_G$, call it $P = u_xw_1\cdots w_ku_y$, as well as any leaves (in $T_G$) of any vertex of $P$ (including $u_y$).

Choose a set $S$ according to the following rule: for all end support vertices $s \in G$, excepting the end support vertex which corresponds to $r$, choose a leaf $\ell_s$ of $s$ in $T_G$ along with the two neighbors in $G$ of $\ell_s$ that are on $C$. Color all vertices of $S$ blue. We claim that $S$ is a zero forcing set of $G$. 

Let $B_1$ be a type-1 branch as described above where $xy \in E(T)$ and $P= u_xw_1\cdots w_ku_y$ is the $u_x$-$u_y$ path in $T_G$. Note that $u_x$ is the only white neighbor of $\ell_{u_x}$ in $G$ and therefore $u_x$ will become blue. Furthermore, all leaves of $u_x$ in $T_G$ will eventually turn blue as the leaves of $u_x$ are ordered in a way that they lie on a path contained in the exterior cycle $C$ of $G$. In fact, we may assume that the leaves of $u_x$ in $T_G$ are indexed as $\ell_i\ell_{i+1}\dots \ell_j$ such that this path is on $C$. Also note that $xy \in E(T)$ if and only if $u_xu_y \in E(G)$ or $\deg_H(w_i) = 2$ for $i \in [k]$. Now once all leaves of $u_x$ in $T_G$ are blue, $w_1$ turns blue as it is the only white neighbor of $u_x$. If $w_1$ has degree $2$ in $G$, then $w_2$ will also turn blue. So assume that $w_1$ is adjacent to a leaf of $T_G$. Without loss of generality, we may assume that either all leaves of $w_1$ in $T_G$ lie on a path of the form $\ell_{i'}\ell_{i'+1}\dots\ell_{i-1}$ where $i' < i$ and $\ell_{i-1}$ is adjacent to $\ell_i$ on $C$, or all leaves of $w_1$ lie on two paths, one of the form $\ell_{i'}\ell_{i'+1}\dots\ell_{i-1}$ where $i' < i$ and $\ell_{i-1}$ is adjacent to $\ell_i$ on $C$ and the other of the form $\ell_{j+1}\dots \ell_t$ where $\ell_j$ is adjacent to $\ell_{j+1}$ on $C$. In either case, one can easily verify that since all vertices in $N_G[u_x]$ are blue, all leaves of $w_1$ in $T_G$ will eventually turn blue at which point $w_2$ will also turn blue. This propagation will continue through $B_1$ and indeed all vertices on a type-1 branch will eventually turn blue. 

Next, choose a type-2 branch $B_2$ as described above where $xy \in E(T)$, neither of $x$ or $y$ is a leaf of $T$, $P = u_xw_1\dots w_ku_y$ is the $u_x$-$u_y$ path in $T_G$ and $x$ has maximum distance from $r$ among all such type-2 branches. Note that this implies that $u_x$ is blue as it is necessarily on a type-1 branch. Since we have assumed $x$ has maximum distance from $r$, we know that $u_x$ is contained in only one type-2 branch. Therefore, if $u_x$ is not adjacent to any leaves in $T_G$, then $w_1$ is the only white neighbor of $u_x$ and propagation continues throughout $B_2$. Therefore, we shall assume that $u_x$ is adjacent to leaves in $T_G$. As above, we may assume that either all leaves of $u_x$ in $T_G$ lie on a path $\ell_{i}\dots \ell_j$ on $C$ or all leaves of $u_x$ lie on two paths, namely $\ell_i\dots \ell_j$ and $\ell_m\dots \ell_n$ on $C$. Choose a type-1 branch $B_1$ that contains $u_x$ as well as a vertex $t$ such that $t$ has a leaf $\ell_{j+1}$ where $\ell_{j+1}$ is adjacent to $\ell_j$ on $C$. Since all leaves of $t$ in $T_G$ are assumed to be blue, all leaves of $u_x$ of the form $\ell_i\dots \ell_j$ will eventually turn blue as well. Similarly, all leaves of $u_x$ of the form $\ell_m\dots \ell_n$ will eventually turn blue. Thus, $w_1$ will be the only white neighbor of $u_x$ at some point and will eventually turn blue as well. Propagation will continue throughout $B_2$ and eventually all vertices of $B_2$ will turn blue. Continuing this same argument by choosing type-2 branches $B_2$ where $xy \in E(T)$ and $x$ has the maximum distance from $r$ among all vertices corresponding to a type-2 branch containing white vertices, one can see that all vertices of $G$ will indeed turn blue and $S$ is a zero forcing set of $G$.

\end{proof}

There are two immediate consequences to the above result. 

\begin{cor} For any Halin graph $G$ where $ES(G)$ is the set of end support vertices in $G$,  $\ggr(G) \geq n(G) - 3(|ES(G)| - 1)$.
\end{cor}

\begin{cor} If $G$ is Halin graph which contains no vertex $v$ with $\deg'(v) \geq 3$, then $Z(G) = Z_{\dot{\ell}}(G) = 3$ and $\ggr(G) = n(G)-3$.
\end{cor}

\begin{proof} Note that in such a Halin graph, there are only two end support vertices and therefore $Z_{\dot{\ell}}(G) \le 3$. On the other hand, $\delta(G) \ge 3$ and by Proposition~\ref{prop:mindegree} $\ggr(G) \le n(G) -3$ and therefore $Z_{\dot{\ell}}(G) \ge 3$. 
\end{proof}

%%%%%%%%%%%%%%%%%%%%%%%%%%%%%%%%%%%%%%%%%%%%%%%%%%%%%%%%%%%%%%%%%
\section{Claw-free cubic graphs}\label{sec:cubic}
As mentioned earlier, Davila and Henning proved the following regarding the zero forcing number of claw-free cubic graphs. 

\begin{thm}\cite{zerocubic}\label{thm:zerocubic}  If $G$ is a claw-free cubic graph with $n(G) \ge 10$, then $Z(G) \le \frac{1}{3}n(G) +1$. 
\end{thm}
We originally wanted to study the loop zero forcing number of claw-free cubic graphs using the same approach to that used in \cite{pdcubic}. As it turns out, we only ever considered zero forcing sets and therefore we provide an improved bound to the one given in Theorem~\ref{thm:zerocubic} in the case when $G$ is a $2$-edge-connected claw-free cubic graph.  To do so, we need the following known results regarding $2$-edge-connected cubic graphs. Recall that the graph obtained from $K_4$ by removing any edge is referred to as a \emph{diamond}. In what follows, a \emph{string of diamonds} is defined to be a maximal sequence $D_1, \dots, D_k$ of diamonds in which, for each $i \in [k]$, $D_i$ has a vertex adjacent to a vertex in $D_{i+1}$ (modulo $k$) and has exactly two vertices which have degree two. Given a string of diamonds $D_1, \dots, D_k$, we label the vertices of $D_i$ as $u_i, v_i, w_i$, and $y_i$ where $u_i$ is adjacent to a vertex in $D_{i-1}$  and $y_i$ is adjacent to a vertex in $D_{i+1}$. Moreover, we refer to $v_i$ and $w_i$ as the \emph{interior vertices of} $D_i$. The following results are necessary for how we approach zero forcing in $2$-edge-connected claw-free cubic graphs. 

\begin{thm}\cite{Oum} A graph $G$ is $2$-edge-connected claw-free cubic if and only if either
\begin{enumerate}
\item[(i)] $G \cong K_4$,
\item[(ii)] $G$ is a ring of diamonds, or
\item[(iii)] $G$ can be built from a $2$-edge-connected cubic multigraph $H$ by replacing some edges of $H$ with strings of diamonds and replacing each vertex of $H$ with a triangle.
\end{enumerate}
\end{thm}

Recall that a $2$-factor of a graph $G$ is a $2$-regular spanning subgraph of $G$. 

\begin{thm}\cite{Plesnik}\label{thm:multigraph-1factor} Let $G$ be a $2$-edge-connected cubic multigraph where $|V(G)|$ is even. For any arbitrary edge $e \in E(G)$, $G$ contains a $2$-factor containing $e$. 
\end{thm}

 In what follows, if we assume that $\mathcal{C} = C_1 \cup \cdots  \cup C_k$ is a $2$-factor of $G$, we will say that $C_i$ is adjacent to $C_j$ if there exists some $v$ on $C_i$ and some $w$ on $C_j$ such that $vw \in E(G)$. We are now ready to prove Theorem~\ref{thm:cubic}, restated here for ease of reference. 
 \vskip3mm

\noindent \textbf{Theorem~\ref{thm:cubic}} \emph{ Let $G$ be a $2$-edge-connected, claw-free cubic graph. If $G$ is not a ring of diamonds, then $Z_{\dot{\ell}}(G)\le Z(G) \le \left\lceil\frac{5n(G)}{18}\right\rceil+1$. If $G$ is a ring of diamonds, then $Z_{\dot{\ell}}(G) =Z(G) =  \frac{n(G)}{4}+2$. 
}
\vskip3mm

\begin{proof} Note that the statement is clearly true when $G = K_4$. Therefore, we may assume that $n(G) \ge 6$. Suppose first that $G$ is a ring of $k$ diamonds where we enumerate the diamonds as $D_1, \dots , D_k$ such that $D_i$ is adjacent to $D_{i-1}$ and $D_{i+1}$ for $i \in [k]$ and the vertices of $D_i$ are labeled $u_i, v_i, w_i, y_i$ where $u_i$ is adjacent to $v_i$, $w_i$, and a vertex of $D_{i-1}$ and $y_i$ is adjacent to $v_i$, $w_i$, and a vertex of $D_{i+1}$. One can easily verify that $S = \{v_1, y_1, u_2, v_2\} \cup \{v_i:3 \le i \le k\}$ is a zero forcing set of $G$ of cardinality $\frac{n(G)}{4}+2$. To see that this is best possible, suppose there exists a loop zero forcing set of $G$, call it $S$, with cardinality $\frac{n(G)}{4}+1$. $S \cap V(D_i) \ne \emptyset$ for each $i \in [k]$ for otherwise $v_i$ and $w_i$ will never become blue. Reindexing if necessary, we may assume $|S \cap V(D_1)| =2$ and $|S \cap V(D_i)| =1$ for $2 \le i \le k$. One can easily verify that one of $v_1$ or $w_1$ is in $S$. Without loss of generality, assume $v_1 \in S$. If $\{v_1, w_1\} \subset S$. then $u_2 \in S$ and $y_1$ will turn blue as all of its neighbors are in $S$. However, $v_2$ and $w_2$ will never turn blue no matter how many times we apply the color change rule so this case cannot occur. Therefore, either $\{v_1, u_1\} \subset S$ or $\{v_1, y_1\} \subset S$. Without loss of generality, we may assume $\{v_1, y_1\} \subset S$. Thus, $u_2 \in S$ for otherwise the color change rule can never be applied. However, we reach a similar contradiction as before in that $v_2$ and $w_2$ will never become blue. Thus, no such loop zero forcing set $S$ of $G$ exists and we may conclude $Z_{\dot{\ell}}(G) =Z(G) =  \frac{n(G)}{4}+2$. 

Hence, we shall assume $G$ can be built from a $2$-edge-connected cubic multigraph $H$ by replacing some edges of $H$ with strings of diamonds and replacing each vertex of $H$ with a triangle. By Theorem~\ref{thm:multigraph-1factor}, $H$ contains a $2$-factor $\mathcal{C'} = C_1' \cup \cdots \cup C_k'$. Let $G'$ be the graph obtained from $H$ by replacing each vertex of $H$ with a triangle. Note that $G'$ may not be $G$. We can create a $2$-factor $\mathcal{C} = C_1 \cup \cdots \cup C_k$ for $G'$ from $\mathcal{C'}$ by replacing each vertex $u$ of $C_i'$ with a $P_3$ $u_1u_2u_3$ where $u_1u_3$ is also an edge in $G'$ (see Figure~\ref{fig:cubexample}).  We may enumerate the vertices of $C_i$ as $C_i = x^i_1\dots x^i_{n_i}$ such that $3 \mid n_i$ and the subgraph of $G'$ induced by $\{x^i_j, x^i_{j+1}, x^i_{j-1}\}$ for all $j \equiv 1\pmod{3}$ is a triangle. Moreover, there exists a perfect matching $M$  in $G'$ between the vertices of 
\[\bigcup_{i\in[k]} \{x^i_1, x^i_4, x^i_7, \dots, x^i_{n_i-2}\}.\]

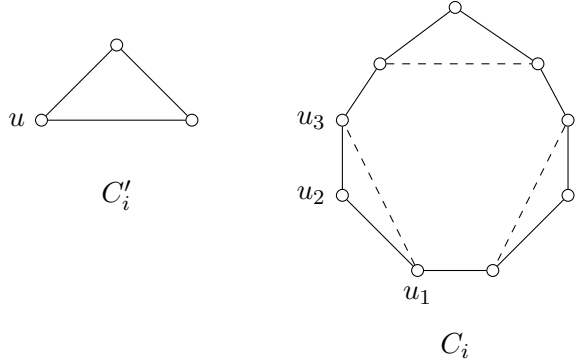
\begin{figure}[h]
\begin{center}
\begin{tikzpicture}[]
\tikzstyle{vertex}=[circle, draw, inner sep=0pt, minimum size=6pt]
\tikzset{vertexStyle/.append style={rectangle}}
	\vertex (1) at (0,0) [scale=.75, label=left:$u$] {};
	\vertex (2) at (2, 0) [ scale=.75] {};
	\vertex (3) at (1, 1) [ scale=.75] {};
	
	\vertex (4) at (5,-2) [scale=.75, label=below:$u_1$] {};
	\vertex (5) at (4, -1) [scale=.75, label=left:$u_2$] {};
	\vertex (6) at (4, 0) [scale=.75, label=left:$u_3$] {};
	\vertex (7) at (4.5, .75) [scale=.75] {};
	\vertex (8) at (5.5, 1.5) [scale=.75] {};
	\vertex (9) at (6.6, .75) [scale=.75] {};
	\vertex (10) at (7, 0) [scale=.75] {};
	\vertex (11) at (7, -1) [scale=.75] {};
	\vertex (12) at (6, -2) [scale=.75] {};
	\node(A) at (1, -1)[]{$C_i'$};
	\node(B) at (5.5, -3)[]{$C_i$};

	\path
		(1) edge (2)
		(2) edge (3)
		(1) edge (3)
		(4) edge (5)
		(5) edge (6)
		(6) edge (7)
		(7) edge (8)
		(8) edge (9)
		(9) edge (10)
		(10) edge (11)
		(11) edge (12)
		(4) edge (12)
		(4) edge[dashed] (6)
		(7) edge[dashed] (9)
		(10) edge[dashed] (12)

	;
\end{tikzpicture}
\end{center}
\caption{Creating a $2$-factor for $G'$ from $\mathcal{C'}$ where the dashed edges are edges in $G'$}
\label{fig:cubexample}
\end{figure}

Let $J$ be the graph with $V(J) = \{u_1, \dots, u_k\}$ where $u_iu_j \in E(J)$ if and only if there exists an edge $xy \in E(G)$ such that $x$ is on $C_i$ and $y$ is on $C_j$, $i\ne j$. Choose a spanning tree $T$ of $J$ rooted at $u_1$  and let $A_i = \{v \in V(T): d_T(u_1, v)=i\}$. Moreover, if $k \ge 3$, choose $T$ and $u_1$ such that $\deg_T(u_1) \ge2$. Let $\ell: V(T) \to [n(T)]$ be any bijection where $\ell(u_1) = 1$ and $\ell(v) < \ell(w)$ if $d_T(u_1, v)<d_T(u_1, w)$. We refer to $\ell(v)$ as the label of $v$ in $T$. Note that we may reindex the vertices of $T$ and corresponding cycles of $\mathcal{C}$ such that $\ell(u_i) = i$ and $u_i$ corresponds to cycle $C_i$ in $\mathcal{C} = C_1 \cup \cdots \cup C_k$.

 Let $M' = \{w_{\alpha_1}v_{\alpha_1}, \dots, w_{\alpha_s}v_{\alpha_s}\}$ be a subset of $M$ where $s = |E(T)|$ and for each $u_{\ell}u_{\ell'} \in E(T)$, there exists $w_{\alpha_j}v_{\alpha_j} \in M'$ such that $w_{\alpha_j}$ is on $C_{\ell}$ and   $v_{\alpha_j}$ is on $C_{\ell'}$. Moreover, we can interchange $w_{\alpha_j}$ and $v_{\alpha_j}$ so that $d_T(u_1, u_{\ell})\le d_T(u_1, u_{\ell'})$. For each $e \in M$, write $e=x^i_{\ell}x_s^j$ where $i \le j$, and if $i=j$, then $\ell < s$.  
 
 Note that since $T$ is rooted at $u_1$, we may choose an edge $e = x_{\ell}^1x_s^2 \in M'$ where $\ell$ and $s$ are as small as possible and reindex the vertices of $C_1$ such that $x_{\ell}^1 = x_1^1$. Similarly, for each $e = x_{\ell}^jx_s^t \in M'$ where $j <t$, we can reindex the vertices of $C_t$ such that $x_s^t = x_1^t$ (and this reindexing will occur only once since $T$ is a tree). Furthermore, if $u_i$ is not a leaf in $T$ and we let  $\{x_{\beta_1}^i, \dots, x_{\beta_j}^i\}$ be those vertices on $C_i$ that are incident to a vertex of the form $x_1^t$ where $t>i$, then we choose the minimum value among $\{\beta_1, \dots, \beta_j\}$ and denote said index by $\eta_i$. Similarly, if $u_i$ is a leaf in $T$, then we let $\{x_{\beta_1}^i, \dots, x_{\beta_j}^i\}$ be those vertices on $C_i$ that are incident to a vertex of the form $x_s^t$ where $t \ne i$. Note that $|\{x_{\beta_1}^i, \dots, x_{\beta_j}^i\}| \ge 2$ as $G$ is $2$-edge-connected. Choose the minimum value strictly greater than $1$ among $\{\beta_1, \dots, \beta_j\}$ and denote said index by $\eta_i$.

 Suppose first that $|V(C_i)| = 6$ for all $i \in [k]$ and let $S = \{x_1^1, x_1^2\} \cup \{x_2^i: i \in [k]\}$. We claim that $S$ is a zero forcing set of $G'$. Note that initially the only white neighbor of $x_1^1$ is $x_6^1$ and therefore $x_6^1$ will become blue. $x_2^1$ will force $x_3^1$ to become blue and $x_6^1$ will force $x_5^1$ to become blue, which will in turn force $x_4^1$ to become blue. Similarly, all vertices of $C_2$ will eventually become blue. Now assume that all vertices of $C_i$ are blue by some timestep $\beta$ where $u_i \in A_{\ell}$ for some $i \in [k]$ and $\ell \in [|E(T)|]$. Let $C_j$ be a cycle in $\mathcal{C}$ where $u_j \in A_{\ell+1}$ and $u_iu_j \in E(T)$. Therefore, there exists $e = x_4^ix_1^j \in M'$. Since $x_4^i$ is blue and $x_1^j$ is its only white neighbor at timestep $\beta$, then $x_1^j$ will be turn blue. Furthermore, $x_2^j \in S$ so it is blue and the propagation will continue through $C_j$. Therefore, $S$ is indeed a zero forcing set of $G'$ of cardinality $\frac{n(G')}{6} + 2$.  Next, note that we can obtain $G$ from $G'$ by replacing some of the edges in $G'$ which were originally in $H$ with strings of diamonds. See Figure~\ref{fig:GfromG'} for an example of how to replace an edge of $G'$ with a string of $2$ diamonds. Choose an interior vertex from each diamond in every string of diamonds used to replace edges in $G'$ to obtain $G$ and call the resulting set $X$. One can easily verify that $S' = S\cup X$ is a zero forcing set of $G$ of cardinality $\frac{n(G')}{6} + 2 + \frac{n(G) - n(G')}{4}$. If $n(G') \ge 12$, then 
 \[|S'| = \frac{n(G)}{4}+2 - \frac{n(G')}{12} \le \frac{n(G)}{4} + 1 \le \frac{5n(G)}{18}+1.\]
 Therefore, we shall assume $n(G') = 6$. If  $n(G)=6$, then $S = N_G(u)$ for any $u \in V(G)$ is a zero forcing set of $G$ of cardinality $3 = \left\lceil\frac{5n(G)}{18}\right\rceil + 1$. So we shall assume that $G$ is obtained from $G'$ by replacing some edges originally in $H$ with strings of diamonds. Suppose first that $G$ is obtained from $G'$ by replacing edges in $G$ which were originally in $H$ with at least three diamonds (i.e. $n(G) \ge 18$). In this case, we label all vertices of each of the $\ell$ diamonds $D_1, \dots, D_{\ell}$ in $G$  as $u_i, v_i, w_i, y_i$ as described above for $i \in [\ell]$. If each edge in $G'$ that was originally in $H$ is replaced with exactly one diamond, then $G$ is isomorphic to the graph depicted in Figure~\ref{fig:small6a} (d). One can easily verify that the black vertices are indeed a zero forcing set. Thus, we shall assume that an edge in $G'$ that was originally in $H$ is replaced with at least two diamonds. Without loss of generality, we may assume $D_1$ and $D_2$ are consecutive diamonds in $G$ where $y_1u_2 \in E(G)$ and $u_1x_5^1 \in E(G)$. One can easily verify that $\{v_i: i \in \ell\}\cup \{y_1, u_2, x_4^1\}$ is a zero forcing set of $G$ of cardinality $3 + \frac{n(G) - 6}{4} \le \frac{5n(G)}{18}+1$ as $n(G) \ge 18$. Therefore, we only need to consider when $n(G) \in \{10, 14\}$.  If $n(G) = 10$, then $G$ is isomorphic to the graph depicted in Figure~\ref{fig:small6a} (a) and the black vertices form a zero forcing set of cardinality $4 = \left\lceil\frac{50}{18}\right\rceil + 1$. If $n(G) = 14$, then $G$ is isomorphic to the graph depicted in Figure~\ref{fig:small6a} (b) or (c). In either case, the black vertices form a zero forcing set of cardinality $5 = \left\lceil\frac{70}{18}\right\rceil + 1$. 
  
  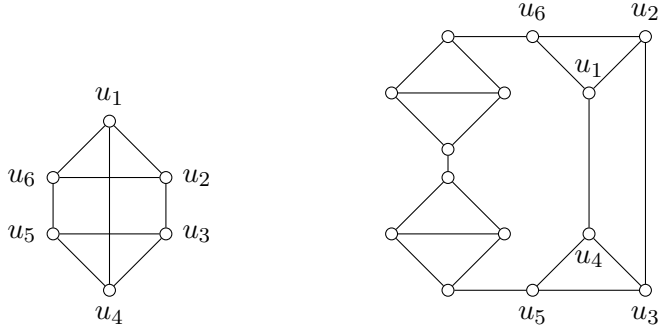
\begin{figure}[h]
\begin{center}
\begin{tikzpicture}[scale=.75]
\tikzstyle{vertex}=[circle, draw, inner sep=0pt, minimum size=6pt]
\tikzset{vertexStyle/.append style={rectangle}}
	\vertex (1) at (0,0) [scale=.75, label=below:$u_4$] {};
	\vertex (2) at (-1, 1) [ scale=.75, label=left:$u_5$] {};
	\vertex (3) at (-1, 2) [ scale=.75, label=left:$u_6$] {};
	\vertex (4) at (0,3) [scale=.75, label=above:$u_1$] {};
	\vertex (5) at (1,2) [scale=.75, label=right:$u_2$] {};
	\vertex (6) at (1,1) [scale=.75, label=right:$u_3$] {};
	
	\vertex (11) at (6, 0) [scale=.75] {};
	\vertex (12) at (5, 1) [scale=.75] {};
	\vertex (13) at (7, 1) [scale=.75] {};
	\vertex (14) at (6, 2) [scale=.75] {};
	\vertex (15) at (6, 2.5) [scale=.75] {};
	\vertex (16) at (5, 3.5) [scale=.75] {};
	\vertex (17) at (7, 3.5) [scale=.75] {};
	\vertex (18) at (6, 4.5) [scale=.75] {};
	\vertex (19) at (7.5, 4.5) [scale=.75, label=above:$u_6$] {};
	\vertex (20) at (9.5, 4.5) [scale=.75, label=above:$u_2$] {};
	\vertex (21) at (8.5, 3.5) [scale=.75, label=above:$u_1$] {};
	\vertex (22) at (8.5, 1) [scale=.75, label=below:$u_4$] {};
	\vertex (23) at (9.5, 0) [scale=.75, label=below:$u_3$] {};
	\vertex (24) at (7.5, 0) [scale=.75, label=below:$u_5$] {};

	\path
	(1) edge (2)
	(2) edge (3)
	(3) edge (4)
	(4) edge (5)
	(5) edge (6)
	(1) edge (6)
	(1) edge (4)
	(2) edge (6)
	(3) edge (5)
		(11) edge (12) 
	(11) edge (13)
	(12) edge (13)
	(12) edge (14)
	(13) edge (14)
	(14) edge (15)
	(15) edge (16)
	(15) edge (17)
	(16) edge (17)
	(16) edge (18)
	(17) edge (18)
	(18) edge (19)
	(19) edge (20)
	(20) edge (21)
	(21) edge (22)
	(22) edge (23)
	(23) edge (24)
	(24) edge (11)
	(19) edge (21)
	(22) edge (24)
	(20) edge (23)

	;
\end{tikzpicture}
\end{center}
\caption{Creating $G$ from $G'$ by replacing an edge with a string of two diamonds}
\label{fig:GfromG'}
\end{figure}

  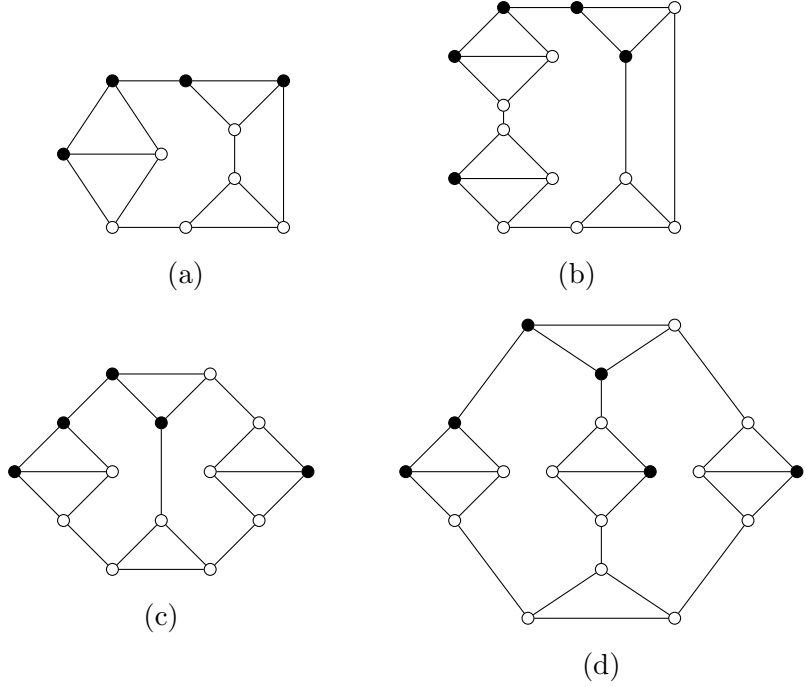
\begin{figure}[h!]
\begin{center}
\begin{tikzpicture}[scale=.65]
\tikzstyle{vertex}=[circle, draw, inner sep=0pt, minimum size=6pt]
\tikzset{vertexStyle/.append style={rectangle}}
	\vertex (1) at (0,0) [scale=.75] {};
	\vertex (2) at (-1, 1.5) [ scale=.75, fill=black] {};
	\vertex (3) at (1, 1.5) [ scale=.75] {};
	\vertex (4) at (0, 3) [scale=.75, fill=black] {};
	\vertex (5) at (1.5, 3) [scale=.75, fill=black] {};
	\vertex (6) at (3.5, 3) [scale=.75, fill=black] {};
	\vertex (7) at (2.5, 2) [scale=.75] {};
	\vertex (8) at (2.5, 1) [scale=.75] {};
	\vertex (9) at (3.5, 0) [scale=.75] {};
	\vertex (10) at (1.5, 0) [scale=.75] {};

	\vertex (11) at (8, 0) [scale=.75] {};
	\vertex (12) at (7, 1) [scale=.75, fill=black] {};
	\vertex (13) at (9, 1) [scale=.75] {};
	\vertex (14) at (8, 2) [scale=.75] {};
	\vertex (15) at (8, 2.5) [scale=.75] {};
	\vertex (16) at (7, 3.5) [scale=.75, fill=black] {};
	\vertex (17) at (9, 3.5) [scale=.75] {};
	\vertex (18) at (8, 4.5) [scale=.75, fill=black] {};
	\vertex (19) at (9.5, 4.5) [scale=.75, fill=black] {};
	\vertex (20) at (11.5, 4.5) [scale=.75] {};
	\vertex (21) at (10.5, 3.5) [scale=.75, fill=black] {};
	\vertex (22) at (10.5, 1) [scale=.75] {};
	\vertex (23) at (11.5, 0) [scale=.75] {};
	\vertex (24) at (9.5, 0) [scale=.75] {};
	
	\vertex (25) at (0, -7) [scale=.75] {};
	\vertex (26) at (-1, -6) [scale=.75] {};
	\vertex (27) at (-2, -5) [scale=.75, fill=black] {};
	\vertex (28) at (0, -5) [scale=.75] {};
	\vertex (29) at (-1, -4) [scale=.75, fill=black] {};
	\vertex (30) at (0, -3) [scale=.75, fill=black] {};
	\vertex (31) at (2, -3) [scale=.75] {};
	\vertex (32) at (1, -4) [scale=.75, fill=black] {};
	\vertex (33) at (3, -4) [scale=.75] {};
	\vertex (34) at (2, -5) [scale=.75] {};
	\vertex (35) at (4, -5) [scale=.75, fill=black] {};
	\vertex (36) at (3, -6) [scale=.75] {};
	\vertex (37) at (2, -7) [scale=.75] {};
	\vertex (38) at (1, -6) [scale=.75] {};
	
	\vertex (39) at (8.5, -8) [scale=.75] {};
	\vertex (40) at (10, -7) [scale=.75] {};
	\vertex (41) at (11.5, -8) [scale=.75] {};
	\vertex (42) at (7, -6) [scale=.75] {};
	\vertex (43) at (6, -5) [scale=.75,fill=black] {};
	\vertex (44) at (7, -4) [scale=.75,fill=black] {};
	\vertex (45) at (8, -5) [scale=.75] {};
	\vertex (46) at (8.5, -2) [scale=.75,fill=black] {};
	\vertex (47) at (10, -3) [scale=.75,fill=black] {};
	\vertex (48) at (11.5, -2) [scale=.75] {};
	\vertex (49) at (10, -4) [scale=.75] {};
	\vertex (50) at (11, -5) [scale=.75,fill=black] {};
	\vertex (51) at (10, -6) [scale=.75] {};
	\vertex (52) at (9, -5) [scale=.75] {};
	\vertex (53) at (13, -6) [scale=.75] {};
	\vertex (54) at (12, -5) [scale=.75] {};
	\vertex (55) at (13, -4) [scale=.75] {};
	\vertex (56) at (14, -5) [scale=.75,fill=black] {};
	
	\node(A) at (1.5, -1)[]{(a)};
	\node(B) at (9.5, -1)[]{(b)};
	\node(C) at (1, -8)[]{(c)};
	\node(D) at (10, -9)[]{(d)};

	\path
	(1) edge (2)
	(1) edge (3)
	(2) edge (3)
	(2) edge (4)
	(3) edge (4)
	(4) edge (5)
	(5) edge (6)
	(6) edge (7)
	(5) edge (7)
	(7) edge (8)
	(8) edge (9)
	(8) edge (10)
	(9) edge (10)
	(10) edge (1)
	(6) edge (9)
	(11) edge (12) 
	(11) edge (13)
	(12) edge (13)
	(12) edge (14)
	(13) edge (14)
	(14) edge (15)
	(15) edge (16)
	(15) edge (17)
	(16) edge (17)
	(16) edge (18)
	(17) edge (18)
	(18) edge (19)
	(19) edge (20)
	(20) edge (21)
	(21) edge (22)
	(22) edge (23)
	(23) edge (24)
	(24) edge (11)
	(19) edge (21)
	(22) edge (24)
	(20) edge (23)
	(25) edge(26)
	(26) edge (27)
	(26) edge (28)
	(27) edge (28)
	(27) edge (29)
	(28) edge (29)
	(29) edge (30)
	(30) edge (31)
	(30) edge (32)
	(31) edge (32)
	(31) edge (33)
	(33) edge (34)
	(33) edge (35)
	(34) edge (35)
	(34) edge (36)
	(35) edge (36)
	(36) edge (37)
	(37) edge (38)
	(38) edge (25)
	(37) edge (25)
	(38) edge (32)
	(39) edge (40)
	(40) edge (41)
	(41) edge (39)
	(39) edge (42)
	(42) edge (43)
	(43) edge (44)
	(44) edge (45)
	(42) edge (45)
	(43) edge (45)
	(44) edge (46)
	(46) edge (47)
	(47) edge (48)
	(46) edge (48)
	(47) edge (49)
	(49) edge (50)
	(50) edge (51)
	(51) edge (52)
	(49) edge (52)
	(50) edge (52)
	(48) edge (55)
	(55) edge (56)
	(56) edge (53)
	(53) edge (54)
	(54) edge (55)
	(54) edge (56)
	(53) edge (41)
	(40) edge (51)

	;
\end{tikzpicture}
\end{center}
\caption{Possible graphs where $\mathcal{C}$ contains only $6$-cycles and $n(G) \le 14$}
\label{fig:small6a}
\end{figure}
  
 Thus, for the remainder of the proof we shall assume that some cycle in $\mathcal{C}$ has length $9$ or more. Suppose first that $\mathcal{C}$ contains at least two cycles. We create a set $D$ as follows. For each $e=x^i_{\ell}x_s^j \in M$ where $i < j$ we place $x_{\ell}^i \in D$ if $|V(C_i)| \ge 9$ and $\ell \ne \eta_i$. For each $e=x^i_{\ell}x_s^j \in M$ where $|V(C_i)| \ge 9$ and  $i=j$, we place  $x_{\ell}^i$ in $D$ if $1 \le \ell < \eta_i$ and $\ell < s$. Otherwise, if $\eta_i < \ell < s$, then we place $x_s^i$ in $D$. We let $S =D\cup \{x_1^2\} \cup \{x_2^i: i \in [k]\}$.  We claim that $S$ is a zero forcing set of $G'$. We proceed by showing that if all vertices of $C_i$ for each $i \in [\alpha]$ where $\alpha \in [k]$ are blue by some timestep $\beta$,  then all vertices of $C_{\alpha+1}$ will eventually be blue. We first show that all vertices of $C_1$ will eventually turn blue. First note that $\{x_1^1, x_2^1, x_1^2, x_2^2\} \subseteq S$. Therefore, all neighbors of $x_1^1$ other than $x_{n_1}^1$ are blue and $x_1^1$ will force  $x_{n_1}^1$ to become blue.  $x_2^1$ will force $x_3^1$ to become blue. Now assume all vertices $x_1^1, x_2^1, x_3^1, \dots , x_j^1$ are blue at some timestep $\beta$ where $j < \eta_1$. If $j \equiv 2\pmod{3}$, then $x_{j+1}^1$ is the only white neighbor of $x_j^1$ and propagation continues in a clockwise fashion. If $j \equiv 0\pmod{3}$, then $x_{j+1}^1 \in S$ of $x_{j+1}^1$ is on some edge $x_{j+1}^1x_s^1 \in M$ where $1 \le s \le j-2$. In either case, $x_{j+1}^1$ is blue or is the only white neighbor of $x_s^1$ which implies $x_{j+1}^1$ will become blue. Finally, if $j\equiv 1\pmod{3}$, then $x_{j+1}^1$ is the only white neighbor of $x_{j-1}^1$ and will therefore become blue. Thus, we may conclude that indeed all vertices of the form $x_j^1$ where $1 \le j \le \eta_1-1$ will eventually be blue. Similarly, moving in the opposite direction from $x_{n_1}^1$ to $x_{\eta_1}^1$, all vertices of the form $x_j^1$ where $\eta_1+1 \le j \le n_1$ will eventually become blue. Hence, either $x_{\eta_1-1}^1$ or $x_{\eta_1+1}^1$ will force $x_{\eta_1}^1$ to turn blue which shows that all vertices of $C_1$ will eventually be blue.

  Now assume that all vertices of $C_i$ for each $i \in [\alpha]$ where $\alpha \in [k]$ are blue by some timestep $\beta$ and consider $C_j$ where $j = \alpha+1$. Therefore, there exists $e = x_r^ix_s^j \in M'$ such that $x_r^i$ is on $C_i$ and $x_s^j$ is on $C_j$ where $i<j$. Moreover, we have indexed $C_j$ such that $x_s^j = x_1^j$. Since $x_r^i$ is blue and $x_1^j$ is its only white neighbor at timestep $\beta$, then $x_1^j$ will become blue. Furthermore, $x_2^j \in S$ so it is blue, $x_{n_j}^j$ will become blue as it is the only white neighbor of $x_1^j$ and $x_2^j$ will force $x_3^j$ to turn blue. If $x_4^j \in S$, then propagation will continue (in a clockwise fashion) around $C_j$. So assume $x_4^j \not\in S$. Thus, $x_4^j$ is on some edge $x_4^jx_s^{\ell} \in M$ where $\ell \le j$. If $j = \ell$, then $x_4^j \in S$, which is a contradiction. Thus, we shall assume that $\ell < j$. By assumption, all vertices of $C_{\ell}$ are blue so $x_4^j$ will also become blue. Hence, propagation will continue until we reach $x^j_{\eta_j}$. Similarly, moving in the opposite direction from $x_{n_j}^j$ to $x_{\eta_j}^j$, all vertices of the form $x_s^1$ where $\eta_j + 1 \le s \le n_j$ will eventually become blue. Hence, either $x_{\eta_j-1}^j$ or $x_{\eta_j+1}^j$ will force $x_{\eta_j}^j$ to become blue which shows that all vertices of $C_j$ will eventually become blue. It follows that $S$ is indeed a zero forcing set of $G'$.

  As noted before,  we can obtain $G$ from $G'$ by replacing some of the edges in $G'$ which were originally in $H$ with strings of diamonds.  Choose an interior vertex from each diamond in every string of diamonds used to replace edges in $G'$ to obtain $G$ and call the resulting set $X$. One can easily verify that $S' = S\cup X$ is a zero forcing set of $G$ of cardinality $|S'| = |D| + k +1+ \frac{n(G) - n(G')}{4}$. We claim that $|S'|\le \left\lceil\frac{5}{18}n(G)\right\rceil+1 $. 
  
  We shall assume for the time being that $k \ge 3$. Assume first that $\mathcal{C}$ contains two cycles $C_i$ and $C_j$ such that $C_i$ and $C_j$ are adjacent and there exist $\ell, \ell' \in [k] - \{i, j\}$ such that $\ell \ne \ell'$, $C_i$ is adjacent to $C_{\ell}$ and $C_j$ is adjacent to $C_{\ell'}$. Moreover, reindexing if necessary, we may assume $i=1$ and $j=2$, $T$ is rooted at $u_1$ (which corresponds to cycle $C_1$), and $\{u_1u_{\ell}, u_2u_{\ell'}, u_1u_2\} \subset E(T)$. Furthermore, we may assume $C_{\ell}$ is chosen so that $x_{\eta_1}^1x_1^{\ell} \in M'$ and $C_{\ell'}$ is chosen so that $x_{\eta_2}^2x_1^{\ell'}\in M'$. Note that each edge in $M$ can be viewed as a matching between triangles in $G'$. Furthermore, for each pair of triangles other than the one containing $\{x_1^1, x_2^1, x_{n_1}^1, x_1^2, x_2^2, x_{n_2}^2\}$, there are at most two vertices from $S'$. Now we can partition $M$ into $M_1$ and $M_2$ such that if $e \in M_1$, then the corresponding triangle pair containing  $e$ contains at most one vertex from $S'$ and if $e \in M_2$, then the corresponding triangle pair containing $e$ contains at least two vertices from $S'$. We claim that $|M_2|-1 \le \frac{n(G')-18}{12}$. Note first that if $e$ is incident to a vertex of the form $x_{\eta_i}^i$, then the corresponding triangle pair contains exactly one vertex from $S'$ and $e \in M_1$. Moreover, if $e \in M_2$, then $e$ is incident to $x_1^s$ for some $s \in [k]$ and the triangle pair associated with $e$ contains vertices of the form $x_j^t, x_{j+1}^t, x_{j-1}^t$ where $t<s$, $j \ne \eta_t$, $x_j^t \in D$ and $x_2^s \in S'$. Enumerate the edges of $M_2$ as $\{e_1, \dots, e_r\}$ where $e_1 = x_1^1x_1^2$. We may assume for $2 \le i \le r$ that  $e_i$ is incident to $x_1^{\alpha_i}$. For each $i \in [r]$, let $T_i$ represent the six vertices of the triangle pair associated with $e_i$. Consider the set $\{W_1, \dots, W_r\}$ where 
  \[W_1 = T_1 \cup \{x_{\eta_1}^1, x_{\eta_1-1}^1, x_{\eta_1+1}^1, x_1^{\ell}, x_2^{\ell}, x_{n_{\ell}}^{\ell}, x_{\eta_2}^2, x_{\eta_2-1}^2, x_{\eta_2+1}^2, x_1^{\ell'}, x_2^{\ell'}, x_{n_{\ell'}}^{\ell'}\}\]
  where $x_{\eta_1}^1x_1^{\ell} \in E(G)$ and $x_{\eta_2}^2x_1^{\ell'} \in E(G)$, and  
  \[W_i = T_i \cup \{x_{\eta_{\alpha_i}}^{\alpha_i}, x_{\eta_{\alpha_i}-1}^{\alpha_i}, x_{\eta_{\alpha_i}+1}^{\alpha_i}, x_1^{\beta_i}, x_2^{\beta_i}, x_{n_{\beta_i}}^{\beta_i}\}\]
  for $2 \le i \le r$ where $x_{\eta_{\alpha_i}}^{\alpha_i}x_1^{\beta_i} \in E(G)$ and $u_{\alpha_i}$ is not a leaf in $T$, and otherwise 
    \[W_i = T_i \cup \{x_{\eta_{\alpha_i}}^{\alpha_i}, x_{\eta_{\alpha_i}-1}^{\alpha_i}, x_{\eta_{\alpha_i}+1}^{\alpha_i}, x_s^{j}, x_{s-1}^{j}, x_{s+1}^{j}\}\]
    for $2 \le i \le r$ where $u_{\alpha_i}$ is a leaf in $T$ and $x_{\eta_{\alpha_i}}^{\alpha_i}x_s^j \in E(G)$ where $u_{\alpha_i}u_j \not\in E(T)$. That is, $x_{\eta_{\alpha_i}}^{\alpha_i}$ is in a triangle pair containing vertices of $C_j$ where the spanning tree $T$ of $J$ does not contain an edge between corresponding vertices $u_{\alpha_i}$ and $u_j$. 
  
  One can easily verify that $W_i \cap W_j = \emptyset$ for $i \ne j$ and by construction, $|W_i| =12$ for $2 \le i \le r$ and $|W_1| = 18$. It follows that $r-1 \le \frac{n(G') - 18}{12}$ which implies $r+2 \le \frac{n(G')}{12} + \frac{3}{2}$. Hence,  $S'$ contains at most one vertex from each triangle pair associated with an edge in $M_1$, there are  $r-1$ triangle pairs that contain exactly two vertices from $S'$, and there is exactly one triangle pair that contains exactly $4$ vertices. Thus, 
  \begin{eqnarray*}
  |S'| &\le& \frac{n(G')}{6} + (r-1) + 3 + \frac{n(G) - n(G')}{4}\\
  &\le& \frac{n(G')}{6} + r+2 + \frac{n(G) - n(G')}{4}\\
  &\le& \frac{n(G')}{6}+\frac{n(G')}{12} + \frac{3}{2} + \frac{n(G) - n(G')}{4}\\
  &=& \frac{n(G)}{4} + \frac{3}{2}\\
  &\le& \frac{5n(G)}{18} +1
  \end{eqnarray*}
  provided that $n(G) \ge 18$. On the other hand, we have assumed $\mathcal{C}$ contains a cycle of length $9$ or more and $k \ge 3$. Thus, $n(G) \ge n(G') \ge 21$. 
  
  Therefore, we shall assume that $\mathcal{C}$ does not contain two cycles $C_i$ and $C_j$ such that $C_i$ is adjacent to $C_j$ and there exist $\ell, \ell' \in [k]- \{i, j\}$ where $C_i$ is adjacent to $C_{\ell}$ and $C_j$ is adjacent to $C_{\ell'}$. It follows that we may write $\mathcal{C} = C_1 \cup \cdots \cup C_k$ where one of the following is true: (a) $k = 3$ and $C_1$ is adjacent to both $C_2$ and $C_3$ and $C_2$ is adjacent to $C_3$, or (b) $k \ge 3$ and $C_i$ is only adjacent to $C_1$ for each $2 \le i \le k$. In case (b), a similar argument to that used above will show $|M_2| -1 \le \frac{n(G')-18}{12}$ where $W_i$ is as defined above for $2 \le i \le r$ and 
  \[W_1 = T_1 \cup \{x_{\eta_1}^1, x_{\eta_1-1}^1, x_{\eta_1+1}^1, x_1^{\ell}, x_2^{\ell}, x_{n_{\ell}}^{\ell}, x_{\eta_2}^2, x_{\eta_2-1}^2, x_{\eta_2+1}^2, x_r^j, x_{r+1}^j, x_{r-1}^j\}\]
  where $x_{\eta_1}^1x_1^{\ell} \in E(G)$ and $x_{\eta_2}^2x_r^j \in E(G)$. Therefore, we assume that we are in case (a). In this case, we may assume that $T = P_3$ where $u_1$ has degree $2$. Thus, $x_{\eta_1}^1x_1^3 \in M'$, $x_1^1x_1^2 \in M'$ and we may assume $x_{\eta_2}^2x_{\eta_3}^3 \in M$. Thus, the triangle pair associated with $x_{\eta_2}^2x_{\eta_3}^3$ does not contain a vertex of $S'$. Moreover, $M_2 = \{e_1\}$ and 
  \begin{eqnarray*}
  |S'| &\le& \frac{n(G')-6}{6} + 3 + \frac{n(G)-n(G')}{4}\\
  &=&\frac{n(G')}{6} + 2 + \frac{n(G) - n(G')}{4}\\
  &=&\frac{n(G)}{4} + 2 - \frac{n(G')}{12}.
  \end{eqnarray*} 
Since we have assumed $\mathcal{C}$ contains a cycle of length $9$ or more and $k=3$, $n(G') \ge 21$ and we have 
\[|S'| \le \frac{n(G)}{4} + 2 - \frac{21}{12} = \frac{n(G)}{4} + \frac{1}{4} \le \frac{5n(G)}{18} + 1.\]
  
  Finally, assume that $ k \in \{1, 2\}$. If $k=2$, then again we may assume $x_{\eta_1}^1x_{\eta_2}^2 \in M$, $M_2 = \{e_1\}$, and the associated triangle pair contains no vertex of $S'$. Moreover, $n(G') \ge 15$ as we have assumed $\mathcal{C}$ contains a cycle of length $9$ or more and $k \ge2$. Thus, 
\begin{eqnarray*}
|S'| &\le& \frac{n(G')-6}{6} + 3 + \frac{n(G) - n(G')}{4}\\
&=& \frac{n(G')}{6} + 2 + \frac{n(G) - n(G')}{4}\\
&=& \frac{n(G)}{4} + 2 - \frac{n(G')}{12}\\
&\le& \frac{n(G)}{4} + 2 - \frac{15}{12}\\
&=& \frac{n(G)}{4} + \frac{3}{4}\\
&\le& \frac{5n(G)}{18} + 1.
\end{eqnarray*}
  Lastly, we shall assume $k=1$ and let \[S = \{x_j^1: j \equiv 1\pmod{3} \text{ and } x_j^1x_r^1 \in M \text{ with } r > j\} \cup \{x_{n_1}^1, x_2^1\}.\]Thus, 
  \begin{eqnarray*}
  |S'| &=& \frac{n(G')}{6} + 2 + \frac{n(G) - n(G')}{4}\\
  &=& \frac{n(G)}{4} + 2 - \frac{n(G')}{12}.
  \end{eqnarray*}
  If $n(G')\ge 12$, then we are done. So assume $n(G') \in \{6, 9\}$. Note that $n(G') \ne 9$ for otherwise there is no perfect matching between $\{x_1^1, x_4^1, x_7^1\}$. On the other hand, if $n(G') = 6$, then $\mathcal{C}$ does not contain a cycle of length $9$ or more, which is a contradiction.

\end{proof}

\end{document}